\documentclass{amsproc}

\usepackage{amssymb,mathrsfs}

\usepackage{xcolor}
\usepackage{tikz}
\usetikzlibrary{decorations.pathreplacing}
\usetikzlibrary{arrows,calc}

\newtheorem{theorem}{Theorem}[section]

\theoremstyle{definition}

\theoremstyle{remark}

\numberwithin{equation}{section}

\begin{document}

\title[Homogenization and Nonlocal Effects]{Multiscale Modeling, Homogenization and Nonlocal Effects: Mathematical and Computational 
Issues}


%    author one information
\author{Qiang Du}
\address{Department of Applied Physics and Applied Mathematics, 
Columbia University, New York, NY 10027,   USA}
\email{qd2125@columbia.edu}
\thanks{The research of Qiang Du is supported in part by the U.S.~NSF grants
DMS-1719699 and CCF-1704833,  AFOSR MURI center for material failure prediction through peridynamics, and ARO MURI Grant W911NF-15-1-0562.
}

\author{Bjorn Engquist}
\address{Department of Mathematics and the Oden Institute, University of Texas, Austin,
TX 78712, USA}
\email{engquist@ices.utexas.edu}
\thanks{The research of Bjorn Engquist is supported in part by the U.S. NSF grant DMS-1620396.}

\author{Xiaochuan Tian}
\address{Department of Mathematics, University of Texas, Austin,
TX 78712, USA}
\email{xtian@math.utexas.edu}
\thanks{The research of Xiaochuan Tian is supported in part by the U.S.~NSF grant DMS-1819233}

\subjclass[2010]{Primary 65-02, Secondary 00A71, 35B27, 64R20, 65N99, 70-08, 74Q99\\
$\mbox{}$\qquad To appear in the {\em Proceedings of the ICERM Symposium Celebrating 75 Years of Mathematics of Computation}, Contemporary Mathematics, AMS
}
%    The 2010 edition of the Mathematics Subject Classification is
%    now available.  If you are citing a classification from the
%    new scheme, use the following input coding instead.
%\subjclass[2010]{Primary }

\date{}

\begin{abstract}
In this work, we review the  connection between the subjects of homogenization and nonlocal modeling
and discuss the relevant computational issues.  By further exploring this connection, we hope to promote the cross fertilization of ideas from the different research fronts. We illustrate how homogenization may help characterizing the nature and the form of nonlocal interactions hypothesized in nonlocal models.  We also offer some perspective on how studies of nonlocality may help the development of more effective numerical methods for homogenization.
\end{abstract}

\maketitle

\tableofcontents

%***********************************Section One********************************************%
%\input intro
\section{Introduction}\label{sec:intro}

Realistic complex systems often involve multiphysics and multiscale processes. Their effective mathematical modeling and simulations are challenging subjects under current study by many scientists. Since the invention of calculus, continuum models in the form of differential equations have been effective models of many physical processes, particularly when the underlying physical quantities are smooth and slowly changing. Mathematical descriptions based on differential operators connect instantaneous rates of changes of relevant quantities
together,  thus are local, since they involve only local information in an  infinitesimal neighborhood. 
Partial differential equations (PDEs) are also used to model rapid and abrupt  changes. An example relevant to our discussion here is given by linear second order elliptic equations with highly-oscillatory and multiscale coefficients. Such equations are popular models of diffusion and transport processes in highly heterogeneous environment, though their analytical studies and numerical simulations are more difficult to carry out, in comparison to equations having smoothly varying coefficients. As possible remedies,  multiscale techniques like homogenization methods have been developed to capture effective properties with either analytical means or computational algorithms based on numerical homogenization \cite{bensoussan2011asymptotic,weinan2003heterognous,weinan2003multiscale,jikov2012homogenization,kornhuber2016numerical,pavliotis2008multiscale}. 
Meanwhile, in recent years, there have been growing interests and capabilities 
in the nonlocal modeling of complex processes that exhibit singularities/anomalies and multiple scales. Examples
include  nonlocal elasticity models like the Eringen model and
peridynamics \cite{BFGS2016,eringen1972nonlocal}. 

The main aim of this work is to explore the connections between the subjects of homogenization and nonlocal modeling
and the relevant computational issues. 
Homogenizations of multiscale problems, both analytical and numerical, are effective ways  to achieve
model reduction and coarse graining, while nonlocality is a generic feature of the latter. We intend to formalize this widely accepted understanding so as to promote the cross fertilization of ideas from the two research fronts. For example,  we hope that  homogenization may help characterizing the nature and form of the nonlocal interactions hypothesized in nonlocal models while studies of nonlocality may help developing more effective numerical methods for homogenization.

Our  discussion is done through a review of some existing works in the literature and some illustrative new examples.
It ranges from the study of deformation of lattice models to
transport, diffusion and wave propagation in highly heterogeneous media, with the goal of building connections between
nonlocal surrogate models and homogenization, model reduction and multiscale modeling. We largely limit ourselves to 
simple linear models to convey the key messages while pointing out some challenges and open questions concerning more complex and nonlinear systems for future research. We start with some general earlier results relating local and nonlocal 
equations, which we have not seen in the homogenization or nonlocal literature. 

\subsection{Local and nonlocal operators: general theorems and dispersion relation}\label{sec:opera}
A differential operator, acting on a function $u$,  is a linear operator that
can be performed locally by taking infinitesimal changes of $u$,
offering its local character. In contrast, nonlocal operators can mean all operators
that are not local.
Let us review some general theorems on linear operators, and 
see that differential operators are the only local operators, and all the other 
linear operators are nonlocal in nature. 
Let $X$ and $Y$ be open sets of $\mathbb{R}^{d_1}$ and $\mathbb{R}^{d_2}$ respectively. 
Any continuous function $k \in C(X\times Y)$ can be viewed as the kernel of an integral operator $L$, 
\[
L u(x) = \int k(x, y) u(y) dy \,.
\]
Conversely, the way to characterize operators having such a kernel is done within the theory of distributions. The {\itshape Schwartz kernel theorem} says that 
there is a  one-to-one correspondence between
distributions $k \in \mathscr{D}^\prime(X\times Y)$ and 
continuous linear maps $L$ from $C_0^\infty(Y)$ and $ \mathscr{D}^\prime(X)$.
$C^\infty_0$ denotes the space of smooth functions with compact support, and
$ \mathscr{D}^\prime$ denotes the space of distributions. 

One may use the support of the kernel $k(x,y)$ to characterize the 
range of interactions of $L$: for $k(x,y)$ that have singular support at $x = y$, $L$ is 
a local operator. For example, $k(x,y)=\delta(y-x)$ leads to
the identity map while $k(x,y)=\nabla \delta(y-x)$
where $\nabla$ denotes distributional gradient operator leads to $L=-\nabla$. On the other hand,  $L$ is a nonlocal operator if the
support of $k$ goes beyond the diagonal. For example, $k(x,y)=\delta(y-x-h)$ leads to a shift operator $Lu(x)=u(x+h)$.
In fact, it was shown in \cite[Chapter 5]{hormander1990analysis} that
the kernel $k$ of a continuous linear map $L: C_0^\infty(X)\to  \mathscr{D}^\prime(X)$ 
is supported by the diagonal if and only if $L$ is in the form of 
 \begin{equation} \label{eq:local-diff}
  L u(x) =\sum a_\alpha(x) \partial^\alpha u(x) \,,
 \end{equation}
 where $a_\alpha \in \mathscr{D}^\prime(X) $ and the sum in \eqref{eq:local-diff} is locally finite. 
 
 An alternative description of an operator $L$ being local is the property that the support of $Lu$ is a subset of the support of $u$, for any suitable function $u$. 
 An interesting discussion, initiated by Peetre \cite{peetre1959} is that differential operators can
 be characterized through such notion of locality. 
 The original {\itshape Peetre theorem} says that a linear operator $L : C_0^\infty(X)\to  C^\infty(X)$ with $\text{supp}(Lu)\subset \text{supp}(u)$ is in the form of \eqref{eq:local-diff} with coefficient $a_\alpha\in C^\infty(X)$. 
 Although the Peetre theorem looks like a variant of the Schwartz's work presented previously, 
 it is actually a deeper theorem. 
First, it says that if the operator $L$ maps to a better space, then the coefficient $a_\alpha$
in \eqref{eq:local-diff} lives in a better space.
Second, no continuity of the map $L$ is assumed in the Peetre theorem,
so the continuity is an automatic fact out of locality. 
So the Peetre theorem characterizes differential operators only through locality. 
 Generalizations of Peetre theorem can be found in \cite{bratteli1986,slovak1988peetre}. 
In particular, the conclusion of locally finite summation in \eqref{eq:local-diff} can actually be replaced with globally finite summation, 
namely the sum in \eqref{eq:local-diff} is for $|\alpha|\leq p$ with a certain fixed number 
$p$ \cite{bratteli1986}.

 From the above discussions, we see that a linear operator is either local or nonlocal.
 Moreover, local operators are differential operators of finite order. 
 Let us also remark that the notions of local and nonlocal in practical modeling are relative to the objects that they operate on: the equation $- \Delta u= v$ is local, but $(-\Delta)^{-1} v=u$ can be seen as nonlocal. Indeed, generically, model reductions of local models often lead to nonlocality and nonlocal models often get localized by closure relations or the introduction of additional state variables \cite{du19cbms}. More discussions on these aspects can be found in later sections.

Based on the Peetre theorem, we have an important observation on the Fourier symbol of $L$, or alternatively the dispersion relation associated with $L$.  For a spatially homogeneous operator $L$, that is, the kernel given by
the Schwartz kernel theorem is translation invariant. 
Peetre's result shows that Fourier symbol of $L$ must be a polynomial of the wave number, if the operator $L$ is local.
Furthermore, once the Fourier symbol of $L$ is not a polynomial of the wave number, 
the operator $L$ with an translation invariant kernel must be nonlocal. 
This simple fact out of the Peetre theorem will be used again and again in section \ref{sec:nonloc}. 

To end this subsection, let us also mention that the generalized Laplacian, or more precisely, the general theory of Dirichlet forms is given in 
the work of Beurling and Deny \cite{beurling1959dirichlet}, 
and it is shown that a Dirichlet form may 
involve local, nonlocal and self-interacting part.  It is then demonstrated by Mosco \cite{mosco1994composite}
 that the limits of Dirichlet forms of local differential operators in general may not be local forms. 
Therefore,  given the generic nonlocal limiting forms of the local problems, one may anticipate that nonlocal limits and related numerical algorithms could play important roles in both theoretical analysis and practical applications.  
In the next, we will demonstrate that nonlocality is a generic feature in the context of homogenization. 
More discussions are followed in section \ref{sec:discre}.

\subsection{Homogenization and multiscale modeling}\label{sec:homog}
Homogenization theory is concerned with the averaged behavior of heterogeneous differential equations with rapidly oscillating coefficients. 
In other words, homogenization aims to obtain the macroscopic or ``homogenized'' or ``effective'' equations from systems with a heterogeneous structure at a microscopic scale. Mathematically, the homogenization problem could be described as follows. 
Let $L^\epsilon$ be a family of differential operators indexed by the small parameter $\epsilon$ and $L^\epsilon$ has coefficients that are oscillatory on the $\epsilon$-scale. Consider the problem 
\begin{equation} \label{eq:diffepsilon}
L^\epsilon u^\epsilon = f^\epsilon\,,
\end{equation}
where $f^\epsilon$ is given. The question is to find a homogeneous effective equation 
\begin{equation} \label{eq:homogenized}
\bar{L} u^0 = \bar{f}
\end{equation}
%this could be described by considering a general equation
%\[
%F_\epsilon (u^\epsilon, x)=0 
%\]
%for a heterogeneous medium 
%with a small parameter dependence $\epsilon$ and the solution $u^\epsilon : \mathbb{R}^d\to\mathbb{R}$.
%The question is to find a homogeneous effective equation 
%\[
%\bar{F}(u, x)=0\,,
%\]
such that $u^\epsilon\to u^0$ in some topology as $\epsilon\to0$. 

We will consider two classical cases of second order differential equations with oscillatory coefficients. The first is one dimensional and the corresponding homogenized equation will be used in section \ref{sec:lions} to derive a homogenized equation, which is nonlocal. The second is a multidimensional elliptic equation, which will be used as a model in the discussion of numerical homogenization in section \ref{sec:projec}. 
Let us illustrate our point by the simple one-dimensional example given also in \cite{engquist2008asymptotic}. 
Consider the problem: 
 \begin{equation}\label{eq:elliptic_1d}
\begin{cases}
- (a_\epsilon(x) u^\epsilon_x)_x = f&\text{in } (0,1) \\
 u^\epsilon(0) =  u^\epsilon(1)=0 &\,,
\end{cases}
\end{equation}
where $a_\epsilon(x)= a(x/\epsilon)>0$ and $a$ is 1-periodic. Notice that 
$$a_\epsilon(x) \rightharpoonup  \int_0^1 a(y) dy =: \langle a \rangle $$
and 
$$a_\epsilon(x)^{-1} \rightharpoonup  \int_0^1 a(y)^{-1} dy : = \langle a^{-1} \rangle \neq \langle a \rangle^{-1} .$$
The solution of \eqref{eq:elliptic_1d} is exactly given by 
\[
u^\epsilon(x)= -\int_0^x  a_\epsilon(\xi)^{-1}\left( \int_0^\xi f(\eta) d\eta  + C^\epsilon  \right) d\xi\,,
\]
where $C^\epsilon$ is determined by the boundary condition at $x=1$ and it is given by
\[
C^\epsilon = -  \int_0^1  a_\epsilon(\xi)^{-1} \int_0^\xi f(\eta) d\eta d\xi /  \int_0^1  a_\epsilon(\xi)^{-1} d\xi\,.
\]
By using the fact that $a_\epsilon(x)^{-1} \rightharpoonup \langle a^{-1} \rangle$, the limit of $u^\epsilon$ as $\epsilon\to0$ is exactly given by 
\begin{equation} \label{eq:homosolu_1d}
u^0(x) =  - \langle a^{-1} \rangle \int_0^x  \left( \int_0^\xi f(\eta) d\eta + C  \right) d\xi\,,
\end{equation}
where $C$ is the constant that gives $u^0(1)=0$. It is now obvious that \eqref{eq:homosolu_1d} satisfies the homogenized equation given by
\[
\begin{cases}
- \bar{a} u^0_{xx} = f, \quad \text{ in } (0,1), \\
 u^0(0) =  u^0(1)=0, 
\end{cases}
\]
where $\bar{a}= \langle a^{-1} \rangle^{-1}$.
This result will be used in section \ref{sec:lions} to derive a nonlocal homogenized problem. 

Consider a general elliptic problem 
\begin{equation} \label{eq:elliptic}
\begin{cases}
- \text{div} \left(A(x/\epsilon)\nabla u^\epsilon \right) =f & \text{in } \Omega\subset \mathbb{R}^d \\
u^\epsilon = 0 & \text{on } \partial \Omega\,,
\end{cases}
\end{equation}
where the matrix $A=(a_{ij})_{1\leq i,j\leq d}$ is assumed to be symmetric, bounded measurable and $Y$-periodic, where $Y=[0,1]^d$ denotes the unit cube. In addition, $A$ satisfies the uniform ellipticity condition. By the classical homogenization theory (\cite{bensoussan2011asymptotic}), the homogenized problem is also an elliptic problem
\begin{equation}\label{eq:elliptic_homogenized}
\begin{cases} 
- \text{div} \left(\bar{A} \nabla u^0 \right) =f & \text{in } \Omega\subset \mathbb{R}^d \\
u^0= 0 & \text{on } \partial \Omega\,,
\end{cases}
\end{equation} 
with coefficient matrix $\bar{A} =(\bar{a}_{ij})_{1\leq i,j\leq d}$ given by
\begin{equation} \label{eq:elliptic_homogenized_A}
\bar{a}_{ij} = \int_Y a_{ij} + a_{ik} \frac{\partial \chi_j}{\partial y_k} dy\,,
\end{equation}
where $\chi_j$ is obtained by solving the cell problem
\begin{equation} \label{eq:cell_problem}
\begin{cases}
-\text{div} \left( A(y) \nabla \chi_j \right) = \text{div} \left( A(y) \nabla y_j \right)  & \text{in } Y\\
\int_Y \chi_j dy =0 \text{ and } \chi_j \text{ is } Y\text{-periodic}\,.
\end{cases}
\end{equation}
This is the classical homogenization case where the original heterogeneous differential equation \eqref{eq:elliptic} is local and the resulting homogenized equation \eqref{eq:elliptic_homogenized} remains local. The process is however nonlocal and involves solving the cell problem \eqref{eq:cell_problem} over a finite domain. 

 {More generally, theoretical tools such as $\Gamma$-, $G$-, and $H$-convergence 
are developed for the analysis of broader settings of homogenization, such as those without assumption of periodicity and
those involving nonlinearities, see \cite{dal2012introduction,jikov2012homogenization,tartar2009general}. }
  { We note that
the nonlinear homogenization theory is less developed. Already in \cite{bensoussan2011asymptotic}
nonlinear homogenization is discussed but the analysis follows the standard linear derivation quite closely, resulting in local PDE effective equations. For more modern examples see, for example, the survey \cite{engquist2008asymptotic}. The homogenized equations  there are however also based on local operators.}

%\subsection{Nonlocal homogenization limit}\label{sec:nlhom}

%***********************************Section Two********************************************%
%\input nonloc

\section{Nonlocal homogenization limit}\label{sec:nonloc}
One key motivation of this work is to  further  explore the connection between nonlocality and homogenization/reduction of complex models. We will present in this and the following sections specific examples for which the homogenized operators,  in contrast to  \eqref{eq:homosolu_1d} and
 \eqref{eq:elliptic_homogenized},  are nonlocal even if the original  operators are local. 
The local limits exhibited in  \eqref{eq:homosolu_1d} and
 \eqref{eq:elliptic_homogenized} may thus be seen as very special cases, although 
  these are also cases that have been extensively studied in numerical homogenization \cite{engquist2008asymptotic}. 
We now present some examples where nonlocality is essential, in the context of homogenization. 
Other examples can be found in, for instance, \cite{bellieud2005homogenization,briane2001fibered,camar2005non,yvonnet2014consistent}.
  
%Thus, the study of such nonlocal features explicitly through nonlocal modeling, analysis and computation \cite{du19cbms} may help development more effective algorithms. 

\subsection{Memory effect through homogenization}  \label{sec:tartar}
Homogenization problems of differential equations often involve coefficients that are highly oscillatory
and only weakly convergent. Nonlocality may thus arise due to the fact that nonlinear functions and weak limit may not commute in general. A classical example that often has been cited is given by Tartar in 
\cite{tartar1989nonlocal} which considered the following equation, 
\begin{equation}
\begin{cases}
\displaystyle
\frac{\partial u^\epsilon}{\partial t} +a_\epsilon(x) u^\epsilon(x,t)=0 \,, \\
u^\epsilon(x,0) = v(x).
\end{cases} 
\end{equation}
Its solution is explicitly given by 
\[
u^\epsilon(x, t) = v(x) e^{-t a_\epsilon(x)}\,. 
\]
The main point in \cite{tartar1989nonlocal} is that if $a_\epsilon(x)$ converges only weakly as $\epsilon\to0$, to some $a_0(x)$,  then the weak limit of $e^{-t a_\epsilon(x)}$ may not be given by $e^{-t a_0(x)}$. Instead, the weak limit $u^0$ 
cannot be generically expressed by $e^{-t b(x)}$ with any function $b(x)$ and 
is given as 
\[
u^0 (x,t) = v(x) \int e^{-t\lambda} d \nu_x(\lambda)\,,
\]
for certain family of probability measures $d \nu_x $. 
With a simple example that $a_\epsilon(x) = a(x/\epsilon)$, where $a(x)$ is given by
\begin{equation} \label{tartar_example}
a(x) = 
\begin{cases}
1 & x\in \bigcup_{k\in \mathbb{Z}} (2k, 2k+1)\\
2 & x\in \bigcup_{k\in \mathbb{Z}} (2k-1, 2k)\,,
\end{cases}
\end{equation}
we see easily $u^0= \frac{1}{2}  v(x)  (e^{-t} + e^{-2t})$, which cannot be expressed by
$v(x) e^{-t b(x)}$ with any function $b(x)$. 
In \cite{tartar1989nonlocal},  Laplace transform in $t$ is used to match the limiting
function $u^0$ with the solution of following nonlocal in time equation 
\begin{equation} 
\begin{cases}
\displaystyle
\frac{\partial u^0}{\partial t} (x,t)+b(x) u^0(x,t)+ \int_0^t K(x, t-s)u^0(x, s)ds=0 \,, \\
u^0(x,0) = v(x)\,,
\end{cases} 
\end{equation}
with a coefficient $b=b(x)$ and a memory kernel $K$. For the above example where $a$ is given by \eqref{tartar_example}, 
the memory kernel is found to be an exponential function and the equation for the limiting function $u^0= \frac{1}{2}  v(x)  (e^{-t} + e^{-2t})$ is given as 
\begin{equation} \label{eq:tartar_homogenized}
\begin{cases}
\displaystyle
\frac{\partial u^0}{\partial t} (x,t)+\frac{3}{2} u^0(x,t)-\frac{1}{4} \int_0^t e^{-\frac{3}{2} (t-s)}u^0(x, s)ds=0\,, \\
u^0(x,0) = v(x)\,.
\end{cases} 
\end{equation}
We note that although the resulting equation for $u^0$ appears essentially nonlocal, it can be localized, not
with $u^0$ directly, rather by defining an additional configuration variable.  For example, with 
$$v^0(x,t)=\int_0^t e^{-\frac{3}{2} (t-s)}u^0(x, s)ds,$$
we get 
\begin{equation} \label{eq:tartar_homogenized_system}
\left\{ 
\begin{aligned}
&\frac{\partial u^0}{\partial t} (x,t)=-\frac{3}{2} u^0(x,t)+\frac{1}{4} v^0(x,t)\,, \\
&\frac{\partial v^0}{\partial t} (x,t)=  u^0(x,t)-\frac{3}{2} v^0(x,t)\,,\\
&u^0(x,0) = v(x)\,,
v^0(x,0)=0\,.
\end{aligned}
\right.
\end{equation}
The introduction of an extended configuration space is one of the popular closure techniques that
has also been used quite effectively in localized nonlocal models, such as in viscoelastiticity. More discussions 
on this can be found later in section \ref{sec:moriz}.

\subsection{Partial differential equation} \label{sec:lions}
We simplify an example given in \cite{bensoussan2011asymptotic} demonstrating that nonlocal operator may exist as the homogenized limit of differential equations. 
\cite{bensoussan2011asymptotic} considers an operator $L_1^\epsilon - \frac{\partial^2}{\partial y^2} L_2^\epsilon$ in the cylinder $O\times \mathbb{R} \subset \mathbb{R}^d \times \mathbb{R}$, where $L_1^\epsilon $ and  $L_2^\epsilon $ are two elliptic operators in $O$ with oscillatory coefficients. The homogenized operator is found through Fourier transform in $y$.
It is said in \cite{bensoussan2011asymptotic} that the Fourier symbol, as $\epsilon\to 0$, is in general not a polynomial, and thus corresponds to a nonlocal operator in the $y$ variable. 
Here we given a simple example reminiscent of  the one in \cite{bensoussan2011asymptotic} and show that nonlocal operator indeed exists as the limit of differential equations.
Consider the following equation
\begin{equation} \label{eq:lions}
- (a_\epsilon(x) u^\epsilon_x)_x + u^\epsilon_{xxyy} = f\quad \text{in } (0,1)\times \mathbb{R}\,,
\end{equation}
where $u^\epsilon : (0,1)\times \mathbb{R} \to  \mathbb{R}$ satisfies zero Dirichlet boundary condition in the $x$ variable and periodic in the $y$ variable. Let  
\[
a_\epsilon = a(x/\epsilon)+1\,,
\]
where $a$ is 1-periodic and 
\[
|a(x)|\leq c_0 <1\,.
\]
Denote $\hat u (k)$ to be the Fourier transform of $u$ in the $y$ variable, 
 then we have
\begin{equation} \label{eq:lions_fourier}
- (a_\epsilon \hat u^\epsilon_x)_x -k^2 \hat u^\epsilon_{xx} = - ((a_\epsilon  + k^2) \hat u^\epsilon_x)_x = \hat f
\end{equation}
Let $b_\epsilon(k) =a_\epsilon  + k^2 $, then use the notations and the 1d homogenization result presented in Section \ref{sec:homog}, the homogenized coefficient $\bar{b}(k)$ for \eqref{eq:lions_fourier} is given by 
\begin{equation} \label{eq:lions_coef_1}
\begin{split}
\bar{b} (k) &=  \langle (a+1 +k^2 )^{-1} \rangle^{-1}  \\
&= \left\langle  (1+k^2)^{-1} \left(1+\frac{a}{1+k^2}\right)^{-1}\right\rangle^{-1}  \\
& = (1+k^2)  \left\langle  1+\sum_{j=1}^\infty \left( \frac{-a}{1+k^2}\right)^j \right\rangle^{-1}  \\
& = (1+k^2)    \left( 1+\sum_{j=1}^\infty (-1)^j \frac{\left\langle a^j\right\rangle}{(1+k^2)^j} \right)^{-1} 
\end{split} 
\end{equation}
Note that if the following equation is true
\begin{equation} \label{eq:lions_aver}
\left\langle a^j\right\rangle = \left\langle a \right\rangle^j \quad \forall j \in \mathbb{N}\,,
\end{equation}
then the right hand side of the last formula in \eqref{eq:lions_coef_1}  
is equal to $1+k^2+\langle a \rangle$. However, \eqref{eq:lions_aver} is generally not true for $j\geq2$.
Indeed, if $a(x)$ is nonnegative, then from H\"{o}lder's inequality, we have always 
$\langle a \rangle^j \leq \left\langle a^j\right\rangle$ for $j\geq2$ and the equality holds only if $a(x)$ is a constant almost everywhere. In general, we have
the following expansion based on the equation \eqref{eq:lions_coef_1}, 
\begin{equation*} %\label{eq:lions_coef_2}
\begin{split}
\bar{b}(k) %&= (1+k^2)    \left( 1+\sum_{j=1}^\infty (-1)^{j} \frac{\left\langle a^j\right\rangle}{(1+k^2)^j} \right)^{-1} \\
&=   (1+k^2)    \left( 1+ \sum_{n=1}^\infty \left( \sum_{j=1}^\infty (-1)^{j-1} \frac{\left\langle a^j\right\rangle}{(1+k^2)^j} \right)^n\right) \\
&=  (1+k^2)   \left( 1+\frac{\langle a \rangle}{1+k^2} + \sum_{j=2}^\infty\frac{1}{(1+k^2)^j} \sum_{d=1}^j \sum_{\substack {l_1+\cdots+ l_d=j \\ l_i\in \mathbb{N}}} \prod_{i=1}^d (-1)^{l_i-1} \langle a^{l_i}\rangle\right) \\
&=   1 +k^2 + \langle a \rangle +  \left(  \sum_{j=2}^\infty\frac{1}{(1+k^2)^{j-1}} \sum_{d=1}^j (-1)^{j-d}\sum_{\substack {l_1+\cdots+ l_d=j \\ l_i\in \mathbb{N}}} \prod_{i=1}^d  \langle a^{l_i}\rangle\right) \,.
\end{split}
\end{equation*}
The term in the parentheses on the right hand side of the last equation given above can be seen as the nonzero correction if \eqref{eq:lions_aver} is violated. In general, $\bar{b}(k)$, as expressed above, is not a polynomial of $k$.
Therefore, by the discussions in section \ref{sec:opera},  the homogenized equation of \eqref{eq:lions} contains a nonlocal operator in the $y$ variable. 

\subsection{Nonlocal effective wave equation}\label{sec:wavdisp}
So far we have been discussing homogenization in a classical sense, namely that 
we exemplify the process of finding the limiting function $u^0$ of $u^\epsilon$ as $\epsilon\to0$
as well as the equation it satisfies. 
From a broader view of homogenization, it is not necessary that the effective equation of 
a highly heterogeneous equation given as \eqref{eq:diffepsilon} should be an $\epsilon$-independent equation 
given as \eqref{eq:homogenized}.
It satisfies the purpose of homogenization as long as one could find an approximation of \eqref{eq:diffepsilon} 
that does not carry the oscillatory behavior and thus is easier for numerical implementation. 
In this spirit, we give in the following an example where an $\epsilon$-dependent effective equation 
 is derived for wave equation in heterogeneous media, and in fact it serves as a better approximation
 of the original problem than the one 
 given by the classical homogenization theory. 
Consider wave propagation through a periodic medium, which is given by the
the equation 
\begin{equation} \label{eq:wave}
\partial_t^2 u^\epsilon(x,t) =  \text{div} \left(A(x/\epsilon)\nabla u^\epsilon \right) \quad x\in \mathbb{R}^d \,,
\end{equation}
complemented with initial condition
\begin{equation} \label{eq:wave_init}
u^\epsilon(x, 0) = f(x),\quad \partial_t u^\epsilon(x,0)=0\,.
\end{equation}
The coefficient matrix $A$ satisfies the same assumptions in section \ref{sec:homog}. 
 Application of the homogenization result in section  \ref{sec:homog} to the spatial part 
 gives the following effective wave equation
\begin{equation}\label{eq:wave_homogenized}
\partial_t^2 u^0(x,t) = \text{div}  (\bar{A}  \nabla u^0(x,t ))\,,
\end{equation} 
where $\bar{A}$ is given by \eqref{eq:elliptic_homogenized_A}. 
It turns out that the effective model \eqref{eq:wave_homogenized} gives a good approximation of \eqref{eq:wave} for short times of observation \cite{bensoussan2011asymptotic,brahim1992correctors,dohnal2014,santosa1991}. 
In a recent work \cite{lin2019uniform}, the authors showed that
for initial-Dirichlet boundary value problems,  $\| u^\epsilon(\cdot, t) - u^0(\cdot,t)\|_{L^2}= O(\epsilon)$
 for $t$ in a fixed time window, independent of $\epsilon$. % \in [0,T]$ with some fixed $T$.  
However, when the size of the time window is large, the important dispersive feature of \eqref{eq:wave} is not 
predicted by the the effective model \eqref{eq:wave_homogenized}. Santosa and Symes \cite{santosa1991} were the first that gave a dispersive 
 effective model by using Bloch wave expansion.
Loosely speaking,  their model is of the form  $\partial_t^2 u =\bar{L}^\epsilon u$, 
with $\bar{L}^\epsilon \approx \Delta + \epsilon^2 (\Delta)^2$, 
and the effective model can well approximate \eqref{eq:wave} when time scale $t\sim O(\epsilon^{-2})$.  
 %See also \cite{allaire2016comparison,conca2002bloch} for more discussions on homogenization via  Bloch wave analysis. 
  There are two drawbacks of such an effective model. 
  The first is that the equation $\partial_t^2 u =\Delta + \epsilon^2 (\Delta)^2$ is actually ill-posed due to the $(\Delta)^2$ term, although 
  it can be made well-posed by a classical Boussinesq trick (\cite{christov1996well,fish2002nonlocal}). The second is that although the dispersive effective model in \cite{santosa1991}
  does better than the non-dispersive effective model for large time scale, the approximation is still within the time scale $O(\epsilon^{-2})$. 
  
  Our main point here is that by allowing the time $t\to\infty$ for given $\epsilon>0$, it is necessary to have 
  a spatially homogeneous  nonlocal operator $\bar{L}^\epsilon$ 
such that the equation
\begin{equation} \label{eq:wave_nonlocal}
\partial_t^2 \bar{u}^\epsilon(x,t)  = \bar{L}^\epsilon  \bar{u}^\epsilon (x,t) \quad x\in \mathbb{R}^d
\end{equation}
complemented with the initial condition \eqref{eq:wave_init} 
 is an approximation of \eqref{eq:wave} for all time.
 Indeed,   $-\bar{L}^\epsilon$ should be the spatially nonlocal operator
 with Fourier symbol corresponding to the first eigenvalue 
 of the operator $-\text{div} \left(A \left(x/\epsilon\right) \nabla u^\epsilon (x,t )\right)$.
And the equation \eqref{eq:wave_nonlocal} is naturally well-posed. 

Let us review the Bloch wave analysis for \eqref{eq:wave} and see how to get an effective equation \eqref{eq:wave_nonlocal}.
Consider the eigenvalue problem for $-\text{div} \left(A \left(x\right) \nabla \right)$: 
\begin{equation} \label{eq:eig_prob}
-\text{div} \left(A \left(x\right) \nabla u(x)\right)  = \lambda u \,.
\end{equation}
For any real vector $k\in \mathbb{R}^d$, there exists a countable number of solutions of \eqref{eq:eig_prob} in the form
\[
\psi_m(x,k)=  e^{2\pi i k\cdot x} \phi_m(x, k) \quad \phi_m(x,k) \text{ is } Y\text{-periodic in } x \,,
\]
with eigenvalue $\lambda_m(k)\in \mathbb{R}$. 
$\psi_m$ is called the quasi-periodic Bloch wave and $0\leq \lambda_0(k)\leq \cdots \lambda_m(k) \leq \cdots \to \infty$. 
The quantity $\sqrt{\lambda_m(k)}$, as a function of k, can be thought of as dispersion relations for the $m$-th mode Bloch wave. 
The Bloch waves form a basis for $L^2(\mathbb{R}^d)$, and
each (complex-valued) function $g \in L^2(\mathbb{R}^d)$ can be expanded 
\[
g(x) = \sum_{m=0}^\infty \int_{Z} \hat g_m (k) \psi_m(x,k) dk,\quad \hat g_m (k)  = \int_\mathbb{R^d} g(x) \bar{\psi}_m(x, k) dx \,,
\]
where $Z=(-\pi,\pi)^d$. Moreover, the Parseval's identity also holds. Please see \cite{bensoussan2011asymptotic} for more details. 
Now the spectral analysis for $-\text{div} \left(A \left(x/\epsilon\right) \nabla \right)$ can be similarly defined by rescaled quantities:
\begin{equation} \label{eig_rescaled}
\psi_m^\epsilon(x,k):= \psi_m(x/\epsilon,\epsilon k), \quad \lambda^\epsilon_m(k) := \frac{1}{\epsilon^2} \lambda_m(\epsilon k)\,.
\end{equation}
In the end, the solution $u^\epsilon$ to \eqref{eq:wave} is written in the expansion \cite{santosa1991}:
\[
u^\epsilon (x,t)= \sum_{m=0}^\infty \int_{Z/\epsilon} \hat f_m (k) \psi^\epsilon_m(x,k) e^{\pm i t \sqrt{\lambda^\epsilon_m(k)}} dk\,.
\] 
There are two key observations in \cite{santosa1991}. 
The first observation is that eigenmodes with $m\geq 1$ can be neglected, in the sense that by defining 
\begin{equation}
u^\epsilon_0 (x,t)= \int_{Z/\epsilon} \hat f_0 (k) \psi^\epsilon_0(x,k) e^{\pm i t \sqrt{\lambda^\epsilon_0(k)}} dk\,,
\end{equation}
we have $\| u^\epsilon(\cdot, t) - u^\epsilon_0(\cdot, t)\| = O(\epsilon) $ for $t\in (0,\infty)$.
The second observation is that $u^\epsilon_0$ can be further simplified by replacing $ \hat f_0 (k)$
and $ \psi^\epsilon_0$  with the 
Fourier transform $\hat f(k)$ and the Fourier mode $e^{2\pi ik\cdot x }$ respectively. 
Let 
\begin{equation} \label{eq:wave_nonloc_solu}
\bar{u}^\epsilon (x,t)= \int_{Z/\epsilon} \hat f (k) e^{2\pi ik\cdot x} e^{\pm i t \sqrt{\lambda^\epsilon_0(k)}} dk\,,
\end{equation}
then we have $\| u^\epsilon_0(\cdot, t) - \bar{u}^\epsilon(\cdot, t) \| = O(\epsilon)$ for $t\in (0,\infty)$.
See \cite{dohnal2014} for more discussions on the spatial norm in these estimates.
Now if the initial condition \eqref{eq:wave_init} is chosen such that  $\hat f (k)$ is supported on $K \subset Z/\epsilon$,
then it is obvious that \eqref{eq:wave_nonloc_solu} is the solution to \eqref{eq:wave_nonlocal} with initial condition \eqref{eq:wave_init}, 
where $- \bar{L}^\epsilon$ is an operator with Fourier symbol that matches $\lambda^\epsilon_0(k)$ for $k\in Z/\epsilon$. 
More specifically, if we let 
\[
 \bar{L}^\epsilon u =\int \gamma^\epsilon(y-x) (u(y)-u(x)) dy\,,
\]
then the kernel $\gamma^\epsilon$ should be chosen such that for $k\in Z/\epsilon$, the equality holds
\begin{equation}\label{eq:wave_nonloc_eig} 
\int \gamma^\epsilon(s) (1-e^{2\pi i k\cdot s }) ds = \lambda^\epsilon_0(k)\,.
\end{equation}
In fact, in order for \eqref{eq:wave_nonloc_eig} to be satisfied, one only need to determine a kernel $\gamma$ according to 
\begin{equation}\label{eq:disp}
\int \gamma(s) (1-e^{2\pi i k\cdot s }) ds = \lambda_0(k)\quad k\in Z\,,
\end{equation}
then by the scaling in \eqref{eig_rescaled}, the kernel $\gamma^\epsilon$ can be determined by a rescaling of $\gamma$:
\[
\gamma^\epsilon(s)=\frac{1}{\epsilon^{d+2}} \gamma\left(\frac{s}{\epsilon}\right)\,.
\] 
Since $\lambda_0(k)$ is always real, we can further assume that $\gamma(s)$ is an even function. With out loss of generality,
if we normalize the integral of $\gamma(s)$ to be $1$, and extend $\lambda_0(k)$ smoothly to $1$ for $k$ large outside $Z$,  then
we may find the function  $\gamma$ from its Fourier transform determined by \eqref{eq:disp}. The smoothness of $1- \lambda_0(k)$
for large $k$ naturally indicates the possibility of a kernel $\gamma$ with a compact support or fast decay. 
Hence, the equation
 \eqref{eq:disp}  in fact gives us a way to construct nonlocal models out of the
dispersion relation of heterogeneous materials and there could be more than one solution to \eqref{eq:disp}. 

Notice that, in the case discussed here, we have
\begin{equation}
\label{eq:Leps}
 \bar{L}^\epsilon u = \frac{1}{\epsilon^{d+2}} \int   \gamma\left(\frac{y-x}{\epsilon}\right)  (u(y)-u(x)) dy\, .
\end{equation}

From the discussions in section \ref{sec:opera}, in order for $\bar{L}^\epsilon$ to be a local operator, it is necessary that
its Fourier symbol is a polynomial of wave number $k$. However, 
$\lambda^\epsilon_0(k)$ in general is not a polynomial of $k$, 
for example, \cite{lee1973waves} shows the dispersion relation for 1d composite materials 
where it is expressed in terms of trigonometric functions. 
Therefore $ \bar{L}^\epsilon$ is in general a nonlocal operator. 
In \cite{conca2002bloch}, $\lambda^\epsilon_0(k)$ is written as a Taylor expansion around $k=0$ 
where the coefficients involving solving a sequence of cell problems. 
It is also found in \cite{conca2002bloch} that all the odd powers in the expansion of $\lambda^\epsilon_0(k)$ vanish.
One can also extend $\lambda_m(k)$ to the complex plane, and it 
is shown in \cite{kohn1959analytic} that for the 1d Schr\"{o}dinger operator (which is similar in our case for 1d), $\lambda_m(k)$, as functions of complex variable $k$, are branches 
of multivalued analytic functions. 
General queries for the analytic properties of periodic elliptic operators can be found in \cite{kuchment2016overview}. 

In addition, let us remark that
%In fact, if one assumes $\gamma$ is integrable with $c$ being the integral of  $ \gamma$, then by smoothly extending the value of $\lambda_0(k)$
%outside $Z$ to match $c$ asymptotically,  one could find the corresponding $\gamma$ by the inverse Fourier transform. 
%And through this way the
%nature of nonlocality could be characterized, namely whether the nonlocal interaction has fast decay and/or with a characteristic length-scale/horizon, once $\lambda^\epsilon_0(k)$ is known.
%Furthermore, 
one may study robust discretization
schemes for \eqref{eq:wave_nonlocal} with the nonlocal operator of the form \eqref{eq:Leps} in the spirit of asymptotically compatible schemes with respect to the 
parameter $\epsilon$  \cite{td14sinum}.

 {
At last, 
it is interesting to see some striking similarities between the analysis above regarding wave propagation operator and the derivation of absorbing boundary conditions in \cite{engquist1977absorbing} even if homogenization is not involved in the latter. Absorbing or far field or radiation boundary conditions are used to artificially limit the size of a computational domain without changing the solution of a PDE. 
In  \cite{engquist1977absorbing}, such accurate conditions were derived in terms of pseudodifferential operators, which are nonlocal in the space-time boundary.  The same concept is applicable to nonlocal problems 
 \cite{du2018numerical} where the notion of boundary conditions are appropriated extended to accommodate nonlocal interactions \cite{du19cbms,du12sirev}.
For PDEs, a hierarchy of local approximative boundary conditions based on differential operators can be given for computational efficiency, which are of increasing order for more accurate approximations of the nonlocal operator  \cite{engquist1977absorbing}. Similar to the discussion above, some of these local approximations give ill posed problems and some well posed ones. For absorbing boundary conditions higher order Taylor expansion of the relevant dispersion relation may give ill posed initial-boundary value problems but Pad\'{e} approximations could lead to well posed problems.}

%***********************************Section Three********************************************%
\section{Nonlocality through  model reduction, numerical homogenization}

Nonlocality is a generic feature of model reduction \cite{du19cbms}. Mathematically this can be easily seen from simple examples like boundary integral formulation of elliptic equations.  This is also evident in many applications, for example, the effective properties, being mechanical, thermal or electrical,  of a heterogeneous medium may be strongly influenced by nonlocal effect \cite{camar2005non,mosco1994composite}.

We will give a couple of examples from numerical homogenization and coarse graining of linear models to illustrate how nonlocality arises naturally and further characterize the nature of the nonlocal interactions.

\subsection{Nonlocality through projection, Schur complements}\label{sec:projec}

We discuss some recent development of numerical homogenization for which projection and Schur complements play an important role \cite{engquist2002wavelet,engquist2008asymptotic,gallistl2017computation, maalqvist2014localization}.
The projection-based numerical homogenization starts with a finite dimensional  
approximation of the differential equation of \eqref{eq:homogenized} on a very fine grid, and tries to project the solution onto the a coarser scale. 
In \cite{dorobantu1998wavelet,engquist2002wavelet},  the coarse scale system is found by using the Schur complement of the fine scale system. It is known that the Schur complement of a sparse matrix, obtained from the finite dimensional approximation of differential operators, is no longer a sparse matrix. So the coarse scale matrix is essentially a nonlocal operator at the the discrete level. 
On the other hand, it is also observed in \cite{dorobantu1998wavelet} that the coefficients of the resulting dense matrix has exponential decay measured by the distance to the diagonal. 
More recently, \cite{maalqvist2014localization}  uses a generalized finite element method based on orthogonal subspace decomposition and the method is re-interpreted in \cite{gallistl2017computation} involving a discrete integral operator.  
In the basic version of of their method, the basis functions have a global support but decay exponentially so it motivates the use of localized basis functions and therefore their method is called the  localized orthogonal decomposition (LOD) method.
Clearly, the two groups of studies both point to the fact that numerical homogenization results in a discrete nonlocal operator with exponential decay in the involved integral kernel under certain assumptions. Here we will establish the fact that the basic version of the subspace decomposition method in \cite{maalqvist2014localization} is in fact exactly the same as the projection method in \cite{dorobantu1998wavelet} using Schur complements.  
Note that the discussions in these works are also related to the variational multiscale method \cite{hughes1998variational}. 

Consider a variational problem defined on a fine space $V$, namely $u\in V$ with
\[
a(u, v) =(f, v) \quad \forall v\in V\,,
\]
where $a$ is a symmetric bilinear form. 
It  is also written as
\[
L u =f \,,
\]
where $L$ is the corresponding linear functional such that $(Lu, v) =a(u,v)$. 
Now let $V_H$ be a subspace of $V$. Then the decomposition $V= V_H\oplus W$ represents the splitting of a function $V$ into a coarse scale function in $V_H$ and its details in the orthogonal complement $W$.
Following the notations in \cite{dorobantu1998wavelet}, we define $P$ to be
the orthogonal projection from $V$ to $V_H$ and $Q=I-P$. Then the system is written into
\[
\begin{pmatrix} 
A & B \\
C & D 
\end{pmatrix}
\begin{pmatrix} 
Qu \\
Pu
\end{pmatrix}
= 
\begin{pmatrix} 
Qf \\
Pf
\end{pmatrix}\,,
\]
where $A=QLQ, B=QLP, C=PLQ, D=PLP$.
The homogenized operator $\bar{L}$ is the Schur complement given by
\[
\bar{L} = D -CA^{-1}B\,,
\]
and the homogenized right hand side is
\[
\bar{f} = Pf - CA^{-1}Qf\,.
\]
The homogenized equation is given by
\begin{equation} \label{eq:Schur_homogenized}
\bar{L} u_H = \bar{f} \,,
\end{equation}
where $u_H = Pu \in V_H$. 
Now the basic version of the LOD method in \cite{maalqvist2014localization} is as follows. Given any $v_H\in V_H$, the corrector $\mathcal{C} v_H \in W$ is defined such that
\begin{equation} \label{eq:LOD_corrector}
a(v_H- \mathcal{C} v_H, w) = 0 \quad \forall w\in W\,.
\end{equation}
The homogenized problem is now to find $u_H \in V_H$ such that 
\begin{equation} \label{eq:LOD_homogenized}
a((1-\mathcal{C}) u_H, (1-\mathcal{C}) v_H) = (f, (1-\mathcal{C}) v_H)\,.
\end{equation}
Let us make several remarks for \eqref{eq:LOD_homogenized}. 
First, if we define a space 
of functions of the form $(1-\mathcal{C}) u_H$, then \eqref{eq:LOD_homogenized}
can be seen as a generalized finite element method with this special finite dimensional space. 
Second, the corrector function $\mathcal{C} u_H$ is given using a family of global basis functions, and thus \eqref{eq:LOD_homogenized} is interpreted as a discrete integral operator in \cite{gallistl2017computation}.  The process of using local basis functions that  approximates the global ones is called in \cite{maalqvist2014localization,gallistl2017computation} the localization.  Third, using \eqref{eq:LOD_corrector}, \eqref{eq:Schur_homogenized} is equivalent to 
\begin{equation} \label{eq:LOD_homogenized_v2}
a((1-\mathcal{C}) u_H, v_H) = (f, (1-\mathcal{C}) v_H)\,.
\end{equation}

We  now show that \eqref{eq:Schur_homogenized} and \eqref{eq:LOD_homogenized_v2} is identical. In fact, the left and right hand side of \eqref{eq:Schur_homogenized} act on $v_H\in V_H$ gives  
\begin{equation} \label{eq:Shur_homogenized_transform}
\begin{split}
((D-CA^{-1} B) u_H, v_H)& =(Pf-CA^{-1}Qf, v_H)\quad  \Longrightarrow   \\
(PLPu_H, v_H) - (PLQ(QLQ)^{-1}QLP u_H, v_H)&= (f, v_H) - (PLQ(QLQ)^{-1}Q f, v_H) \\
\Longrightarrow \; (L u_H, v_H) - (L (QLQ)^{-1}QL u_H, v_H)&=  (f, v_H) - (L(QLQ)^{-1}Qf, v_H) \,.
\end{split}
\end{equation}
On the other hand, the definition of corrector in \eqref{eq:LOD_corrector} is the same as 
\[
(Lv_H, w) = (L \mathcal{C} v_H, w)\quad  \forall w\in W\,,
\]
and it is equivalent to $QL v_H = QL  \mathcal{C} v_H = QL Q  \mathcal{C} v_H $. 
So  we have
\begin{equation}\label{eq:corrector_expression}
\mathcal{C} v_H =  (QLQ)^{-1} QL v_H \quad \forall v_H\in V_H\,. 
\end{equation}
\eqref{eq:corrector_expression} explains why $\mathcal{C} v_H$ is expressed global basis functions because the inverse of $QLQ$ gives a nonlocal operator. 
Using \eqref{eq:corrector_expression},  the left hand side of \eqref{eq:LOD_homogenized_v2} is given as
\[
\begin{split}
a((1-\mathcal{C})u_H, v_H)&=(L u_H, v_H) - (L \mathcal{C}u_H, v_H) \\ 
&= (L u_H, v_H) - (L(QLQ)^{-1} QL u_H, v_H)\,,
\end{split}
\]
and the right hand side of \eqref{eq:LOD_homogenized_v2} is 
\[
\begin{split}
 (f, (1-\mathcal{C}) v_H)  &= (f, v_H) - (f, (QLQ)^{-1}QLv_H)\\
 &= (f, v_H) - (Qf, (QLQ)^{-1}QLv_H) \\ 
 & =  (f, v_H) - ((QLQ)^{-1} Qf, Lv_H) \\
 &= (f, v_H) - (L (QLQ)^{-1} Qf, v_H)\,,
 \end{split}
\]
and they are the same as the expressions in \eqref{eq:Shur_homogenized_transform}.
This shows that \eqref{eq:Schur_homogenized} is the same as \eqref{eq:LOD_homogenized} and \eqref{eq:LOD_homogenized_v2}. 

 {It is observed in \cite{maalqvist2014localization} that the global basis functions
decay exponentially away from the node they are associated with, which is in accordance with
the observation in \cite{dorobantu1998wavelet} that the non-zero entries in the Schur complement have exponential decay. This justifies the use of local basis functions where they are constructed from cell problems on local patches 
with Dirichlet boundary conditions. This is also related to the oversampling technique in the context of multiscale  
finite element method that will be discussed in the following subsection.
Finally, we remark that the exponential decaying property is true under the assumption that 
the bilinear forms are uniformly bounded and coercive. If the heterogeneous coefficients contain scales with very high contrast, then the exponential decaying property may be lost. This may be due to the fact that the intrinsic nonlocality could be developed under such cases. In fact, 
  \cite{mosco1994composite} gives an example of a sequence of local problems with coefficients that develop singular measures in the limit, and the sequence is shown to converge to a nonlocal problem.  
}

%\subsection{Nonlocaity through projection/Schur complements}\label{sec:projec}

%\input msfem
\subsection{MsFEM and HMM}\label{sec:msfem}

We will briefly discuss the Multiscale Finite Element Method (MsFEM), \cite{efendiev2013generalized,efendiev2009multiscale,hou1997multiscale,
hou2004removing}, and the Heterogeneous Multiscale Method (HMM), \cite{abdulle2016reduced,abdulle2012heterogeneous,weinan2003heterognous, weinan2003multiscale}.

These are two numerical methodologies that also are influenced by analytical homogenization based with cell problems. The methods are however more general and do not explicitly rely on cell problems. The similarity with homogenization is that there is a macroscale $O(1)$ and microscales $O(\epsilon)$ and a nonlocal operation of cell problem type. Most rigorous analysis for convergence of these methods is furthermore done in the setting of heterogeneous differential equations with a microstructure and based homogenization theory as in section \ref{sec:homog}.

\begin{figure}[htbp]
\centering
\includegraphics[width=8cm]{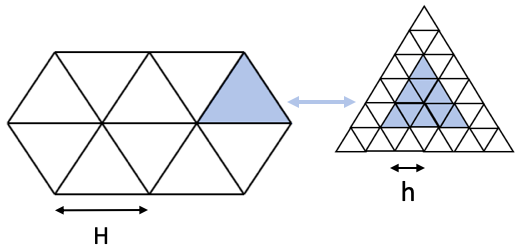}
\caption{MsFEM with the fine scale element size $h$ and the coarse scale element size $H$.}
\label{fig:msfem}
\end{figure}

We will first consider the MsFEM. The goal in this method is to solve a differential equation with a finite element method when the coefficients in the differential equation have a microstructure as in equation \eqref{eq:elliptic}. Standard piecewise polynomial basis functions require element sizes $h$, which are smaller than $\epsilon$ to resolve the oscillations. The key to the success of MsFEM is the development of multiscale basis functions, that are specific for each differential equation. The formulation with these new basis functions allows for efficient approximation even if their element size $H$ is larger than $\epsilon$. The computation of the new basis functions is based on standard finite elements of size $h$ over domains of size $O(H)$ that also approximates the original differential operator. Here $h$ is smaller than $\epsilon$, which is typically smaller than $H$. The relation between $h$ and $H$ is plotted in Figure \ref{fig:msfem}.  {Since it is difficult to determine the accurate boundary conditions for computing the basis functions, the oversampling technique is often used, which means that the $O(H)$ domain above is a bit larger than $H$. The error created by the boundary conditions decays exponentially with the increase of oversampling layers \cite{chung2018fast}.}   For efficiency this step can be done in parallel and ahead of the final finite element solve on the macroscale. The computation of the multiscale elements is a nonlocal operation since a solution operator is involved and thus has strong similarities to the homogenization process as outlined in section \ref{sec:homog}.

\begin{figure}[htbp]
\centering
\includegraphics[width=8cm]{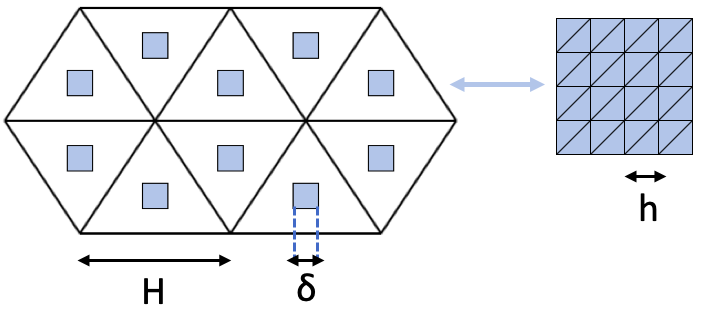}
\caption{HMM with the the fine scale element size $h$, coarse scale element size $H$ and
the size of domain for  mircoscale simulation $\delta$.}
\label{fig:hmm}
\end{figure}

The HMM framework is not restricted to FEM discretization and can, for example, couple a finite volume method on the macroscale and molecular dynamics on the microscale.
For such a case it is similar to the nonlocal version of the quasicontinuum method \cite{tadmor1996quasicontinuum}. We will here for comparison focus on the use of FEM for both micro and macro scale. Finite element HMM (FEM-HMM) is less ambitious than MsFEM in that the target is an approximation of the effective or homogenized solution and not the full multiscale solution. This means computing $u^0$ in \eqref{eq:elliptic_homogenized} instead of $u^\epsilon$ in \eqref{eq:elliptic}. It is however more ambitious in that the microscale simulations are concentrated on small subsets of the full computational domain.
The subsets are of size $\delta$, which is smaller than the coarse scale elements, 
see Figure \ref{fig:hmm}. This means that much more extreme range of scales can be approximated but it also requires much clearer scale separation or gaps in the scales of the full multiscale problem.

Since the target is $u^0$ we will assume that something is known about an effective equation. For example, the structure of \eqref{eq:elliptic_homogenized} is known but not the stiffness matrix $\bar{A}$. It is then possible to formulate a FEM discretization of \eqref{eq:elliptic_homogenized} but for the explicit components of homogenized matrix $\bar{A}$. The size of the elements can be of the order $H$ as in MsFEM. This macroscale simulation is then coupled to a microscale simulation based on \eqref{eq:elliptic} with element size $h$, which provides the missing components that are needed to determine the stiffness matrix. This microscale computation is done on small domains where the components are required typically around the location of the quadrature points in the macroscale FEM. 
 {The micro problem is usually complemented with periodic or Dirichlet boundary conditions. 
A more thoughtful way of imposing boundary conditions for the micro problem is to take into account the far field and create artificial boundary conditions that 
 does not change the solution of a PDE when restricting it to small domains, see e.g. \cite{engquist1977absorbing}. }

\subsection{A linear lattice model as an illustration}
\label{sec:lattice}
Following the discussion on the Schur complement given previously, 
we can further provide a simple illustration using coarse grained models of linear lattice models as done in \cite{duliyuan2019}. More specifically,  \cite{duliyuan2019} considered a lattice model involving  next nearest neighbor interactions associated with the energy,
\begin{equation}
\label{eq: ENext}
E[u] = \sum_{x\in\mathbb{L}} \frac{K_1}{2}\left(\frac{u(x+\varepsilon)-u(x)}{\varepsilon}\right)^2 + \frac{K_2 }{2} \left( \frac{u(x+2\varepsilon) - u(x)}{\varepsilon}\right)^2 - f(x) u(x),
\end{equation}
where  $\mathbb{L}$ denotes a one-dimensional uniform lattice with $\varepsilon$ being the lattice space and $K_1$ and $K_2$ being the force constants satisfying the so-called phonon stability condition
$K_1>0$, and $K_1+4K_2>0$.

The variational equation corresponding to \eqref{eq: ENext} is given by the finite difference equation
\begin{equation}
\label{eq: EqNext}
 -K_1\frac{u(x+\varepsilon)-2u(x)+u(x-\varepsilon)}{\varepsilon^{2} }- K_2 \frac{u(x+2\varepsilon)-2u(x)+u(x-2\varepsilon)}{\varepsilon^{2} } = f(x)
\end{equation}
which can be viewed as a discretization of the Poisson equation
 $- (K_1 + 4K_2) u''(x) = f(x)$ for a one-dimensional elastic bar.
 A coarse lattice can be defined by selecting a representative atom from each group of $M$ atoms, 
 as done in the local version of the quasicontinuum method \cite{tadmor1996quasicontinuum}. Let us denote the selected representative atoms by
\begin{equation}
\label{lcg}
\mathbb{L}_{CG} = \big\{y | y = nH, n\in\mathbb{Z} \big\},\quad  \mathbb{L}^-_{CG} = \mathbb{L} \backslash \mathbb{L}_{CG}
\end{equation}
where $H = M\varepsilon$ represents the spatial scale of interest.
The original model on the full lattice $\mathbb{L}$ can be reduced to a model on $\mathbb{L}_{CG}$ by eliminating
the non-representative atoms in $\mathbb{L}^-_{CG} =\mathbb{L}\setminus \mathbb{L}_{CG}$
that are between neighboring atoms in $\mathbb{L}_{CG}$.
 The reduced model can be written as,
\begin{equation}
\label{eq: reduced}
\sum_{y'\in \mathbb{L}_{CG}} \theta(y-y') u(y') = \bar{f}(y), \;\forall y \in \mathbb{L}_{CG}.
\end{equation}
This procedure is similar to the Schur complement reduction described earlier
in section \ref{sec:projec}. As reflected by properties of the nonlocal 
kernel function $\theta$, a main observation in \cite{duliyuan2019} is
 that the reduced model on the coarse lattice $\mathbb{L}_{CG}$ is nonlocal: the displacement of all the representative atoms are in principle coupled with  every other representative atom. Indeed, it was shown in  \cite{duliyuan2019}
 that the resulting nonlocal interaction through the coarse-graining has no compact support but
 decays faster than any algebraic power, that is, $
\theta(y) = o(|y|^{-s})$ for any $s>0$.  The numerical test in fact suggested that $\theta$ decays exponentially. 
 
 Moreover, in the case where $K_2>0$, 
  the Equation \eqref{eq: EqNext}  preserves the discrete maximum principle so that
   $\theta$ is an even function with a zero sum, positive at the origin but negative away from the it.
Furthermore, with a numerically verified postulation on
$- \sum \theta(y) y^2 = 2M( K_1+4K_2)$, it was shown in \cite{duliyuan2019} that
 there exists a function $\theta_0$ such that,
\begin{equation}
\begin{aligned}
-\frac{1}{2}\sum_{x\in\mathbb{L}} \theta_0(x)x^2 & =
-\frac{1}{2M}\sum_{y\in \mathbb{L}_{CG}} \theta(y)y^2 = -\frac{1}{2}\sum_{x\in \mathbb{L}} M\theta(Mx) x^2,
\end{aligned}
\end{equation}
holds for all sufficiently large $M$, we let $\theta_0(x) = M\theta(M x)$ for $x \in \mathbb{L}$. Then, in
the experiments given in \cite{duliyuan2019}, it was shown that
as $M$ increases, the rescaled function $M\hat\theta(\xi/M)$ indeed converges to a fixed function.
Thus, $M$ serves as a parameter for the characteristic range of the resulting nonlocal interactions derived through
a coarse graining process.

%\subsection{Linear lattice model}\label{sec:lattice}

%\input moriz

%\subsection{Mori-Zwanzig}\label{sec:moriz}

\subsection{Memory effect and nonlocality through  Mori-Zwanzig}\label{sec:moriz}

One can generalize the discussion of model reduction via Schur complement to the dynamic case which naturally
lead to the appearance of nonlocal in time memory effect.
A discussion of the same nature is given by the so-called Mori-Zwanzig formalism which exemplified
the nonlocal memory effect in coarse-grained dynamic systems \cite{mori65,zwanzig01}. The subject is
of much importance in numerical simulations, for example, it is tied to the choice of thermo-stat in molecular dynamics
simulations and dissipative Brownian particle dynamics \cite{junghans2008transport,ceriotti2009langevin}.
More recent discussions and applications can be found in \cite{chorin2000optimal,stinis2015renormalized,li2015incorporation,hijon2010mori}.

Adopting similar notations used before, we let
 $P$ be an orthogonal projection operator 
and $Q=I-P$. Starting from  a simple autonomous system 
$\dot{x}(t)=Lx(t)$ with a linear operator $L$ defined on the state space and
$x(0)=x_0$, 
a direct decomposition leads to a coupled system
$$\frac{d}{dt}\left(
\begin{array}{c}
Px(t)\\
Qx(t)
\end{array}\right)
=\left(
\begin{array}{cc}
PLP & PLQ\\
QLP & QLQ
\end{array}\right)
\left(
\begin{array}{c}
Px(t)\\
Qx(t)
\end{array}\right).
$$
Assuming that the reduced quantity of interest (QoI) is  $u=Px$, 
by elimination and substitution, we can get 
\begin{equation}
\label{eq:reduced}
\dot{u}(t)=PLPu(t)+PLQ\int_0^t e^{(t-s)QLQ}QLPu(s) ds + PLQe^{tQLQ}Qx_0\,.
\end{equation}
The resulting differential integral equation contains
the second term  on the right hand side that
 reflects the nonlocal in time memory effect and 
possible  nonlocal interaction
in the state variables even if $L$ itself is local.
Let us also remark that if $Qx$ is a fast variable, then the nonlocal equation \eqref{eq:reduced} may
be approximated by a differential equation again. 
Take the example in section \ref{sec:tartar}. We know that the reduced equation
of the system \eqref{eq:tartar_homogenized_system} is given by the nonlocal in time equation \eqref{eq:tartar_homogenized}. 
However, if we assume $v^0$ is a fast variable by considering a modified version of \eqref{eq:tartar_homogenized_system}:
\begin{equation} \label{eq:tartar_system_eps}
\left\{ 
\begin{aligned}
&\frac{\partial u^0}{\partial t} (x,t)=-\frac{3}{2} u^0(x,t)+\frac{1}{4} v^0(x,t)\,, \\
&\frac{\partial v^0}{\partial t} (x,t)= \frac{1}{\epsilon} \left(u^0(x,t)-\frac{3}{2} v^0(x,t)\right)\,,\\
&u^0(x,0) = v(x)\,,
v^0(x,0)=0\,,
\end{aligned}
\right.
\end{equation}
then as $\epsilon\to0$, the reduced equation of \eqref{eq:tartar_system_eps} gets back to a local differential equation again. 
This is related to the averaging theorem and more discussions can be found in \cite{pavliotis2008multiscale}.
 
For a nonlinear system, a similar derivation can be made.
We let the coarse grained vairables
 $u(t)$ be the QoI and work with the associated
Liouville equation for $u$
given by $\dot{u}(t)=Lu(t)$ where $L$ is the time-independent
linear Liouville operator.

Formally, we have $u(t)=e^{tL}a$ but this cannot be directly calculated since 
 $L$ contains information on all the state variables, not just $a=u(0)$. 
Thus, we take
$$\dot{u}(t)=Lu(t)=e^{tL}La =e^{tL}
PL a+ e^{tL} QL a,$$
and apply the
Dyson's identity 
$$e^{tL}=e^{t(QL+PL)}=e^{tQL} +  \int_0^t e^{sL}PLe^{(t-s)QL} ds $$
to get
\begin{equation}
\label{eq:mz0}
\dot{u}(t)=e^{tL}
PL a+ \int_0^t e^{(t-s)L}PL e^{sQL}QL a\,ds +
e^{tQL}
QL a.
\end{equation}
This is similar to \eqref{eq:reduced}, with no reduction done yet
from the original system.

Next, by making a choice of $P$, such as the projection operator $Pv=<v, a> <a,a>^{-1} a$ given by Mori \cite{mori65} where
 $<\cdot, \cdot>$ denotes the correlation operator 
(or $P$ can be given as a conditional expectation 
with respect to $a$ \cite{zwanzig01}), we may rewrite the equation \eqref{eq:mz0} as 
\begin{equation}
\label{eq:mz1}
\dot{u}(t)=G u(t) - \int_0^t  \Gamma(s) u(t-s) ds + f(t).
\end{equation}
On the right hand side of \eqref{eq:mz1},
 the first term
$$Gu(t)=e^{tL} PL a =  e^{tL} <La, a><a,a>^{-1} a
=  <La, a><a,a>^{-1}u(t)$$
 gives the Markovian term operating only on the QoI
and depending only on  the current time.
The third term $f(t)=e^{tQL}QLa$ is commonly treated as the stochastic fluctuation. This is where the reduction  begins to take effect. 
The Mori-Zwanzig formalism ensures a thermodynamically consistent approach
through appropriate coupling with the second term.

Specifically, 
by defining a memory kernel $\Gamma(t)$ satisfying 
$$\Gamma(t)a=-< Lf(t), a><a,a>^{-1} a=PLf(t)=
PL e^{sQL}QL a,$$  we see that
 the second term on the right hand side of \eqref{eq:mz1} indeed becomes
$$ -\int_0^t  \Gamma(s) u(t-s) ds =  
- \int_0^t e^{(t-s)L} \Gamma(s) a\, ds
=\int_0^t e^{(t-s)L}PL e^{sQL}QL a \, ds$$
and represents the nonlocal memory, thus non-Markovian, 
 effect in the reduced system.

In the case that $L$ is a skew-symmetric operator (associated with a Hamiltonian system),  we have
$$\Gamma(t) a
 =-< Lf(t), a><a,a>^{-1} a =<f(t), La> 
<a,a>^{-1} a.$$
Since $f(t)$ is in the range of $Q$, this leads to
$$\Gamma(t) a
%= <Q f(t), La><a,a>^{-1} a
= < f(t), QLa><a,a>^{-1} a
=< f(t), f(0)><a,a>^{-1} a.$$

The above equation is  seen as a fluctuation-dissipation relation \cite{kubo1966fluctuation} connecting the memory (dissipation)
effect with the correlation of fluctuation.
Meanwhile, by noticing that $<f(t), a>=0$ due to orthogonality of $P$ and $Q$, we can also get from \eqref{eq:mz1} the relation
$$
\dot{C}(t)=G C(t) - \int_0^t  \Gamma(s) C(t-s) ds$$
 for the time correlation 
$C(t)=<u(t),u(0)>=<u(t), a>$.

From \eqref{eq:mz1}, we see that
the nonlocal nature of the reduced stochastic dynamics is again intrinsic, although localization is often sought after in practice
through approximations of the nonlocal operator.
For example, it is well known that
nonlocal interactions may be localized if the kernel
$\Gamma(t)$ can be approximated by a finite combination
of exponentials.
This can be seen by using the property $\Gamma'(t)=\Lambda \Gamma(t)$, 
which is satisfied by exponential kernels of the type
$\Gamma(t)=e^{\Lambda t} \Gamma_0$ 
(including the Dirac $\delta$-measure at the origin and
its finite order of derivatives in the distribution sense). 
A simple calculation then yields
$$\dot{v}(t)=\Gamma_0 u(t) + 
\int_0^t  \Gamma'(s) u(t-s) ds=
\Gamma_0 u(t) + \Lambda v(t),$$
a formally local in time equation involving the newly introduced variable
$$v(t)=\int_0^t  \Gamma(s) u(t-s) ds.$$
This means that localizaton can be achieved
through the introduction of extended configuration variables,
as demonstrated in the treatment of memory effect in the constitutive modeling 
of viscoelastic materials (via the introduction of Maxwell stress, 
for example).
Other examples of effective localizaton approaches
include methods based on rational approximations and Galerkin projections \cite{li2015incorporation,hijon2010mori,chu2018asymptotic}.

%***********************************Section Four********************************************% 

%\input discre

\section{Nonlocality and homogenization, open questions}\label{sec:discre}

We have discussed several natural connections of nonlocality,  homogenization and model reduction.
This serves to motivate cross fertilization of ideas from the different subjects, particularly for modeling and simulation
of complex phenomena. 
For example, on one hand we discuss
how nonlocal modeling can benefit from works on multiscale model reduction and homogenization
that  helps to infer the nature of nonlocal interactions and the type of nonlocal interaction kernels. 
On the other hand, understanding nonlocal models that could be the  homogenized limit of local PDEs
could provide insight to the design and analysis of numerical homogenization of PDEs.

\subsection{Characterizing nonlocality through  homogenization and model reduction}

While nonlocal models have become popular in recent years, their rigorous derivations have remain limited. 
Often, the nonlocal interactions are postulated and the choices are validated against experiments.

Homogenization and coarse graining, as illustrated in this work, could be valuable tools to help characterizing
the nature of nonlocality. There have been previous discussions along this line, for example, in the 
context of deriving peridynamics from molecular dynamics or other fine scale models \cite{lehoucq2007statistical,si2014origin}, but rigorous mathematical theory remains largely missing. Moreover, while there are extensive studies on how to find conditions to yield local limit for local forms, it is interesting to see 
when to expect the limiting nonlocal forms that have exponential decay or with compact support as illustrated in some earlier examples. These issues may in turn offer insight on the design of robust and effective numerical algorithms for complex multiscale problems.

  {
There are a few cases of nonlinear problems where the rigorously derived effective equations 
from homogenization  and coarse  graining
have nonlocal components. For example,
Tartar considers the nonlinear Carleman equations with highly oscillatory initial data in \cite{tartar1980}. 
There is little rigorous treatment of more realistic physical and engineering models but peridynamics serves as an
 example of a more empirical nonlocal model that captures the nonlinear crack formation and elastic interaction in brittle materials. An alternative can be the quasicontinuum multiscale method mentioned in Sections \ref{sec:msfem},
which is more based on first principles but with prohibitively high computational cost for practical applications.}

\subsection{Homogenization of nonlocal models}
Homogenization techniques for
local continuum PDE counterpart can be useful in the study of nonlocal models such as peridynamic models of
cracks and fractures. 
Meanwhile, we can also develop homogenization theory for nonlocal models.
 A number of works have studied homogenization problems of linear nonlocal models for which the nonlocal interaction
 kernel display multiscale features, see for example,  \cite{alali2012multiscale,du2016multiscale,mengesha2015multiscale}. 
In these works, however, the multiscale features are prescribed on a finer scale than the horizon parameter for the
nonlocal range of interactions,
thus leading to generically nonlocal homogenized limits with the same range of nonlocal horizon parameter.
It could be interesting to examine cases where the multiscale features of the kernels are coupled with the horizon parameter.
Then, it might be possible to end up with homogenized limits with nonlocal interactions of a different nature, or for example, interactions described by a different horizon parameter. In the latter case, it is conceivable that the limit could also be local. 
Hence, the local homogenized limit of nonlocal problems and the nonlocal homogenized limit of local problems may then be seen as special instances of more general classes of homogenization problems. A possibly unified perspective could in turn
 offer further insights to the special cases.

  What could also be studied further is to develop homogenized models of nonlinear nonlocal mechanics: while
  resolving the nonlocal models allows us to find solutions with  large deformations and singularities like cracks, it could be very  
  interesting to see if the nonlocal models can be homogenized or coarse grained to yield
  models that  can capture the averaged mechanical responses in the presence of complex micro cracks
  without the need to track the detailed crack patterns.

\subsection{Exploring nonlocality for robust and asymptotically compatible discretizations} \label{subsec:AC}

\begin{figure}[htb!]
\centering
       \begin{tikzpicture}[scale=0.9]
     \tikzset{to/.style={->,>=stealth',line width=.9pt}}
 \node(v1) at (0,4.2) {\textcolor{red}{\Large $u_{\epsilon}^{h}$}};
  \node (v2) at (5.8,4.2) {\textcolor{red}{\Large $u_{0}^{h}$}};
   \node (v3) at (0,0) {\textcolor{red}{\Large $u_\epsilon$}};
    \node (v4) at (5.8,0) {\textcolor{red}{\Large $u_0$}};
         \node[align=center,yshift=-0.2cm,anchor=east] at (v1.west)
        {Discrete\\Nonlocal};    
           \node[align=center, yshift=0.2cm,anchor=east] at (v3.west)
        {Continuum\\Nonlocal};   
      \node[align=center, yshift=-0.2cm,anchor=west] at (v2.east)
        {Discrete\\Local};      
                   \node[align=center,yshift=0.2cm,anchor=west] at (v4.east)
        {Continuum\\Local PDE}; 
    \draw[to] (v1.east) --  node[midway,below] {\textcolor{blue}{$\epsilon\to 0$}}
     (v2.west);
    \draw[to] (v1.south) -- 
     node[midway,xshift=0.3cm, yshift=0.5cm,rotate=-90, right]  {\textcolor{blue}{$h\to0$}}
       (v3.north);
     \draw[to] (v3.east) --  node[midway,above] {\textcolor{blue}{$\epsilon\to 0$}}
      (v4.west);
    \draw[to] (v2.south) --  node[midway,xshift=-0.3cm, yshift=0.5cm,rotate=-90,right]  {\textcolor{blue}{$h\to0$}}
    (v4.north); 
   \draw[gray,dashed][to] (v1.south east) to[out = -55, in = 175, looseness = 1.2] node[midway,xshift=0.1cm,yshift=0.6cm,rotate=-38] 
    {\textcolor{blue}{$\epsilon\to0$}  } node[midway, xshift=0.5cm, yshift=1.1cm,rotate=-42] 
    {\textcolor{blue}{$h\to0$}  }
  (v4.north west);
    \draw[gray, dashed][to] (v1.south east) to (v4.north west);
     \draw[gray, dashed][to] (v1.south east) to[out = -5, in = 120, looseness = 1.3]   (v4.north west);
      \end{tikzpicture}
      \vskip-10pt
         \caption{A diagram of possible paths between $u_\epsilon$, $u_\epsilon^h$, $u_0^h$ and $u_0$ via various limits.}
   \label{dt-fig:diagram}
\end{figure}
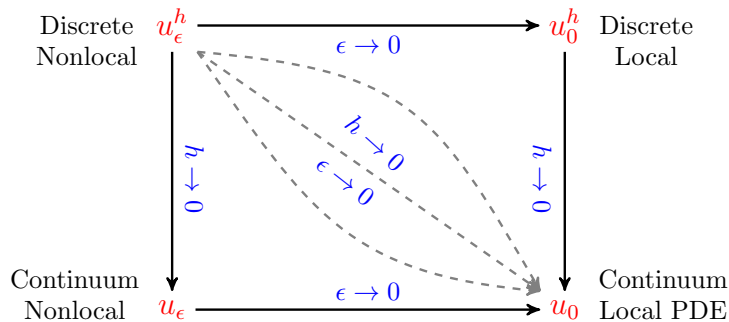

For homogenization problems that are characterized by a small length scale $\epsilon$,
 one might be interested in examining the dependences of numerical approximations of the solution
 on the parameter $\epsilon$ including but not limited to uniform convergence in $\epsilon$.
 Let us take an example in nonlocal modeling, where the parameter $\epsilon$, called the horizon, stands for the 
  length of nonlocal interaction. 
   One is interested in robust approximations
insensitive to model parameter $\epsilon$, and the notion of {\em AC scheme} \cite{td13sinum,td14sinum}  serves
such  a purpose.
More specifically, let $u_\epsilon$ be the solution to the nonlocal problem associated with the horizon parameter
and the $\epsilon\to0$ limit is given by the solution $u_0$ to the corresponding local problem.  Take $h$ as the discretization parameter (mesh spacing or the level of numerical
 resolution). We denote by $u_\epsilon^h$ the discrete approximation of $u_\epsilon$ and
 $u_0^h$ the discrete approximation of $u_0$,  the limit of  $u_\epsilon$ as $\epsilon\to 0$, 
 then the diagram given in Figure \ref{dt-fig:diagram}, as introduced in \cite{td14sinum},
 serves to highlight the different paths connecting the different solutions and limits.
A key ingredient to assure the AC property of approximations to nonlocal variation problems
considered in \cite{td14sinum} is the notion of asymptotic denseness property of the finite dimensional approximation
spaces. 

Figure \ref{dt-fig:diagram} and the research behind it can be generalized as inspiration for future exploration. 
Instead of starting with a nonlocal model and end up with a local one, 
we could start with a local problem with $\epsilon$ representing a fine scale parameter and consider
the limiting process as $\epsilon\to0$. And as we saw in sections \ref{sec:intro} and \ref{sec:nonloc} this process may end with a local problem or a nonlocal problem.  Indeed, the discussion in \cite{camar2005non,mosco1994composite}
implies that nonlocal models might be more generic as limit of certain classes of homogenization problems. 
 {As previously mentioned, local problems with heterogeneous coefficients with scales of very high contrast may be approximating 
a nonlocal problem in the limit.  
 In this case, it is reasonable that  one should take into account of the nonlocal nature of the limit
 and design basis functions that are encoded with information of the limiting nonlocal problem to get an AC scheme.}

Another direction of generalization is to include a numerical scale $\delta$, which is larger than $\epsilon$.
In Figure \ref{dt-fig:diagram} the horizon in the nonlocal formulation serves both as defining the model and
the domain for the numerical computation. If the nonlocal kernel is defined on the whole space without a compact support, then
$\epsilon$ could denote an essential width of the kernel and $\delta$ the numerical horizon that would be chosen larger
for higher accuracy. Indeed, the approximation to the fractional Laplacian using the numerical horizon is discussed in \cite{delia2013fractional,tian2016acm}. 
In MsFEM  $\delta$ would be $O(H)$ and in HMM $\delta<H$. In both cases we will have $h<\epsilon<\delta<O(1)$. 
As in the discussion related to Figure  \ref{dt-fig:diagram} it will be advantageous to have asymptotic 
compatibility with respect to different limits of all these scale variables. 
For homogenization problems of local PDEs with a fine scale parameter,
it is well known
that the careful design of finite dimensional basis are also crucial components of effective numerical homogenization algorithms, see discussions in the review \cite{engquist2008asymptotic}. 
More extensive explorations of such connections remain to be carried out in the future. 
In particular, 
research effort in this regard might help the design of more effective schemes for homogenization without separation of
scales.

\subsection{Reducing the computational complexity}
There are on-going studies on reducing the computational complexity involved in multiscale problems.
For example,
precomputing certain quantities for later use in the overall algorithm has been successfully applied in techniques related to numerical homogenization in order to reduce the overall computational cost. Some examples are precomputation of the multiscale elements in MsFEM \cite{chung2016adaptive,efendiev2013generalized,efendiev2009multiscale}, and computing the effect of the generalized cell problems in HMM \cite{abdulle2012heterogeneous} ahead of the macroscale simulation. The latter is also labeled sequential HMM or parametrization. The efficiency can be further enhanced by relying on reduced basis techniques. See, for example \cite{abdulle2016reduced}.

Since the computational cost per unknown is typically very high in nonlocal modeling it would be interesting to see if some of the techniques mentioned above can be transferred from numerical homogenization to more general nonlocal models. A recent development, somewhat in this direction, involves precomputing, following the methodology of the fast multipole method in lowering the computational complexity in nonlocal diffusion models \cite{tian2019fast}.

 {There have also been some efforts in numerical 
homogenization in the spirit of asymptotic compatibility when $\epsilon$ and $\delta$ 
approach each other. The goal is to reduce the, so called, resonance error from
$O(\epsilon/\delta)$ to $O((\epsilon/\delta)^\gamma)$ for $\gamma>1$ and thereby 
allowing for smaller $\delta$ with the same error tolerance 
\cite{arjmand2016time,gloria2016reduction}. 
Another possibility for lowering the computational complexity is to reduce the nonlocal 
domain size  in both numerical homogenization and nonlocal models. 
This connects to the discussion in section \ref{subsec:AC} regarding the size of the horizon when 
approximating fractional Laplacians. 
In addition, seamless coupling of local and nonlocal models through spatially varying horizon parameters and heterogeneous localization
may also provide effective remedies for the reduction of complexity \cite{du2019seam,tian2017trace}. More studies are needed to utilize such techniques for
adaptive computation.
}

%\section{Nonlocality and discretization, open questions}\label{sec:discre}

\bibliographystyle{amsplain}

\bibliography{ref}

\end{document}